# Visualizing Rotations and Composition of Rotations with Rodrigues' Vector

**Angel G. Valdenebro**
Instituto de Estructura de la Materia (CSIC).
Serrano 123, 28006 Madrid. Spain

**Abstract**
The purpose of this paper is to show that the mathematical treatment of three-dimensional rotations can be simplified, and its geometrical understanding improved, by using the Rodrigues' vector representation. We present a novel geometrical interpretation of the Rodrigues' vector. Based on this interpretation and simple geometrical considerations, we derive Euler-Rodrigues' formula, Cayley's rotation formula, and the composition law for finite rotations. The level of this discussion should be suitable for undergraduate physics or engineering courses where rotations are discussed.

## 1  Introduction

The mathematical treatment of rotations is of interest to scientists, engineers and computer programmers that deal with rigid structures, physical or virtual, from molecules and crystal lattices to aircrafts, and from robotic arms to computer graphics. One of the most interesting parametrizations of three-dimensional rotations, both from the theoretical and from the computational point of view, is the Rodrigues' vector. However, its use is not as extended as it deserves to be. One possible reason is the lack of a clear geometrical interpretation of the Rodrigues' vector and the related formulas. This is precisely the gap that we would like to fill with this article.

The plan is the following. In section 2 we set notation and review the definition of Rodrigues' vector and its most important merits. We also review some well-known rotation formulas, including brief clarifications of their not so well-known history and not so clear naming. In section 3 we provide a novel geometrical interpretation of Rodrigues' vector. In sections 4 and 5 we use this interpretation to derive, based on simple geometrical considerations, several rotation formulas, and the rotation composition law. Finally, in sections 6 and 7, we present some brief applications of Rodrigues' vector and its geometrical interpretation, to infinitesimal rotations and to rigid body kinematics respectively.

## 2  Review of Basic Facts

Let $\mathbf{R}$ be the rotation matrix that relates a certain vector $\mathbf{x}$ in the three-dimensional Euclidean space and the corresponding rotated vector $\mathbf{x}' = \mathbf{R}\mathbf{x}$. A well-known form of $\mathbf{R}$ as a function of the unit vector $\mathbf{n}$ along the rotation axis and of the rotation angle $\theta$, is

$$\mathbf{R}(\mathbf{n},\theta) = (\cos\theta)\mathbf{1} + (\sin\theta)(\mathbf{n}\times) + (1-\cos\theta)\mathbf{n}\mathbf{n}^T, \qquad (1)$$

where $(\mathbf{n}\times)$ and $\mathbf{n}\mathbf{n}^T$ are the matrix operators that transform the vector $\mathbf{x}$ into $\mathbf{n}\times\mathbf{x}$ and $\mathbf{n}(\mathbf{n}\cdot\mathbf{x})$ respectively, and $\mathbf{1}$ is the 3x3 unit matrix. We have supposed $\mathbf{R}$ to be a dextro-rotation, that is, the orientation of the axis and the rotation are related by the right-hand rule (grabbing the rotation axis with the right hand, the extended thumb points along $\mathbf{n}$ and the other fingers indicate the sense of the rotation.) We have adopted the "active" point of view for rotations: they move the vectors and not the coordinate axes.

This formula was (and still is sometimes) erroneously attributed to Onlinde Rodrigues, who published it [1] before vector notation was developed, relating the coordinates of a point before and after rotation (one equation per coordinate) and using the direction cosines of the rotation axis (that are nothing else but the components of $\mathbf{n}$.) Now it is generally recognized [2] that Euler published a similar formula 65 years earlier [3], relating the direction cosines of a point before and after rotation. For that reason, the formula is usually called Euler-Rodrigues' formula. The formula with modern vector notation was introduced by Gibbs [4].

Defining the vector $\mathbf{Q} = \tan(\theta/2)\mathbf{n}$, Eq. (1) can be transformed into

$$\mathbf{R}(\mathbf{Q}) = \mathbf{1} + \frac{2}{1+(\mathbf{Q}\cdot\mathbf{Q})}\left[(\mathbf{Q}\times) + (\mathbf{Q}\times)^2\right]. \qquad (2)$$

The idea is due to Rodrigues [1], and the components of $\mathbf{Q}$ (remember that vector notation was not available yet) are usually called Rodrigues' parameters. Vector form is found in Gibbs' book [4], who gave the cumbersome name "semi-tangent vector of version" to $\mathbf{Q}$, and no credit to Rodrigues for it. Nowadays $\mathbf{Q}$ is called Rodrigues' vector or Gibbs' vector depending on the source. Formula (2) appears in the literature without an eponymic designation. The just mentioned Rodrigues' parameters should not be confused neither with the four Euler-Rodrigues symmetric parameters (or Euler symmetric parameters, or quaternion components) [2,5] nor with the three modified Rodrigues parameters [5].

It may seem that Rodrigues' vector is just a mathematical trick to hide the angle dependence into the norm of the axis vector, and simplify rotation formulas, but it is not so. It has remarkable properties. It provides a parametrization of three-dimensional rotations based on just three independent parameters, while the Euler-Rodrigues formula depends on four parameters (the rotation angle and the three components of $\mathbf{n}$) that are not independent (the norm of $\mathbf{n}$ has to be 1.) Furthermore, it provides a one-to-one mapping between $\mathbb{R}^3$ vectors and rotations of an angle in the open interval $(-\pi,\pi)$. Although Rodrigues' vector $\mathbf{Q}$ is not defined for the $\pi$ value of $\theta$, rotations of an angle $\pi$ can be treated as a limiting case [6]. Rodrigues' vector corresponding to the null rotation (unit rotation matrix) is the zero vector. Rodrigues' vector corresponding to the inverse of a given rotation is the negative of the Rodrigues' vector of the given rotation. In the limit of infinitesimal rotations, Rodrigues' vector becomes infinitesimal, and the related formulas can be simplified by neglecting terms of second (or higher) order in $\mathbf{Q}$. As we will see in the next paragraphs, Rodrigues' vector is closely related to Cayley transform (in three dimensions) and is the key for one of the most interesting formulas for the composition of finite rotations.

Arthur Cayley proved the following result [7]. For any orthogonal matrix $\mathbf{R}$ which does not have -1 as an eigenvalue, there is one, and only one, skew-symmetric matrix $\mathbf{A}$, that allows us to express $\mathbf{R}$ as

$$\mathbf{R} = (\mathbf{1}-\mathbf{A})^{-1}(\mathbf{1}+\mathbf{A}). \tag{3}$$

Reciprocally, the matrix $\mathbf{A}$ can be written as

$$\mathbf{A} = (\mathbf{R}-\mathbf{1})(\mathbf{R}+\mathbf{1})^{-1}. \tag{4}$$

Matrices $\mathbf{A}$ and $\mathbf{R}$ are said to be the Cayley transform of each other. Orthogonal matrices with determinant equal to 1 (special orthogonal matrices) are (isomorphic to) proper rotations. Cayley's rotation formula (3) is important because it is valid in any dimension, and parametrizes *n*-dimensional rotations in terms of $n(n-1)/2$ independent parameters (the components of an *n*-dimensional skew-symmetric matrix.) In three dimensions, $\mathbf{A}$ is precisely the skew-symmetric matrix obtained from the Rodrigues' vector, $\mathbf{A} = (\mathbf{Q}\times)$, that is,

$$\mathbf{R} = (\mathbf{1}-\mathbf{Q}\times)^{-1}(\mathbf{1}+\mathbf{Q}\times), \tag{5}$$

and

$$(\mathbf{Q}\times) = (\mathbf{R}-\mathbf{1})(\mathbf{R}+\mathbf{1})^{-1}. \tag{6}$$

The exclusion of three-dimensional rotations with -1 as eigenvalue is equivalent to excluding rotations of angle $\pi$.

The composition of consecutive rotations is easily expressed using Rodrigues' vectors,

$$\begin{aligned}\mathbf{R}(\mathbf{Q}_3) &= \mathbf{R}(\mathbf{Q}_2)\mathbf{R}(\mathbf{Q}_1), \\ \mathbf{Q}_3 &= \frac{\mathbf{Q}_1 + \mathbf{Q}_2 + \mathbf{Q}_2 \times \mathbf{Q}_1}{1 - \mathbf{Q}_2 \cdot \mathbf{Q}_1}.\end{aligned} \tag{7}$$

If the denominator of the second expression in Eq. (7) vanishes the resulting rotation is through an angle $\pi$ about an axis parallel to $\mathbf{Q}_1 + \mathbf{Q}_2 + \mathbf{Q}_2 \times \mathbf{Q}_1$. The rotation composition law (7) is elegant, simple and computationally economic. Some interesting properties of the composition of rotations can be easily obtained from Eq. (7). Finite rotations do not commute. The order in which the rotations are composed changes the axis of the resultant rotation (the direction of $\mathbf{Q}_3$) but not its angle (the norm of $\mathbf{Q}_3$). To derive these properties, note that $\mathbf{Q}_2 \times \mathbf{Q}_1$ and $\mathbf{Q}_1 \times \mathbf{Q}_2$ are different, opposite to each other, and perpendicular to $\mathbf{Q}_1 + \mathbf{Q}_2$, while $\mathbf{Q}_2 \cdot \mathbf{Q}_1 = \mathbf{Q}_1 \cdot \mathbf{Q}_2$.

An eloquent defense of further merits of Rodrigues' vector, stressing its relevance for crystallography, can be found in Refs. [8] and [9].

The abundance of papers [10–14] and classic books [15] providing geometrical derivations of expression (1) demonstrates how important it is for students to gain a geometrical understanding of rotation formulas. However, there are neither geometrical interpretations of the Rodrigues' vector (other than the definition itself) nor direct geometrical derivations of expressions (2), (5), (6) and (7). As we said in the introduction, this is the gap that we try to fill in the next sections.

## 3  Geometrical Interpretation of Rodrigues' Vector

The geometrical interpretation of Rodrigues' vector, or more precisely, of the skew-symmetric matrix operator obtained from it, is grounded on the following propositions (a), (b) and (c), illustrated on Figs. 1a, 1b and 1c respectively.

Consider a dextro-rotation of angle $\theta$ about an oriented axis through the origin of coordinates. Let $\mathbf{Q}$ be the Rodrigues' vector of this rotation.

(a) The operator $(\mathbf{Q}\times)$ applied to any position vector $\mathbf{x}$ results in a vector $(\mathbf{Q}\times\mathbf{x})$ which is tangent to the arc described by point $\mathbf{x}$ as it rotates about the axis, and connects $\mathbf{x}$ to the bisector of the rotation angle.

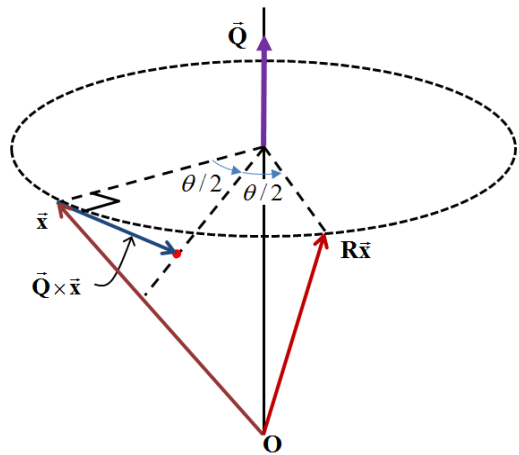

*Fig. 1a. Geometrical interpretation of $(\mathbf{Q}\times)$*

(b) The operator $(\mathbf{1}+\mathbf{Q}\times)$ applied to $\mathbf{x}$ results in the position vector of the intersection point of the above mentioned tangent and bisector.

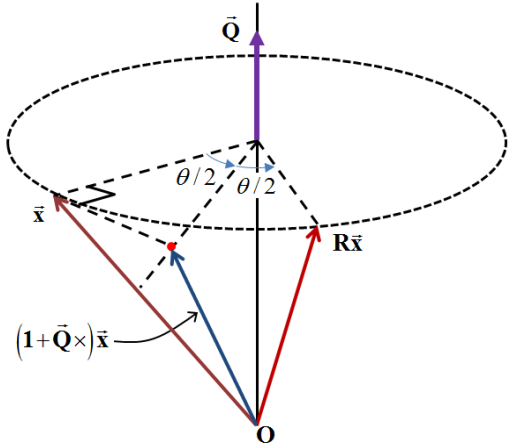

*Fig. 1b. Geometrical interpretation of $(1+\mathbf{Q}\times)$*

(c) For a unit vector **a** located on the plane perpendicular to the rotation axis, the vector $(1+\mathbf{Q}\times)\mathbf{a}/\|(1+\mathbf{Q}\times)\mathbf{a}\|$ corresponds to the point located on the rotation arc, at an angular distance from **a** equal to half the rotation angle.

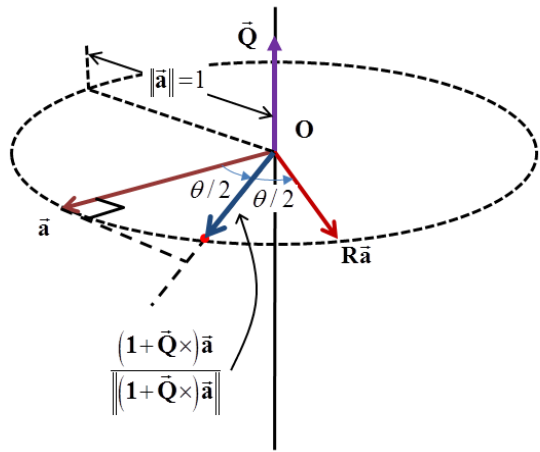

*Fig. 1c. Geometrical interpretation of $(1+\mathbf{Q}\times)\mathbf{a}/\|(1+\mathbf{Q}\times)\mathbf{a}\|$ being **a** a unit vector perpendicular to the rotation axis.*

Proposition (a) is easily justified based on three observations that require only simple trigonometry and knowledge of the cross product properties. The direction of $\mathbf{Q}\times\mathbf{x}$ is tangential to the rotation arc at $\mathbf{x}$, because it has to be perpendicular to $\mathbf{x}$ and to $\mathbf{Q}$ (i.e.: to the rotation axis.) The radius of the rotation arc has length $\|\mathbf{Q}\times\mathbf{x}\|/\|\mathbf{Q}\|$ (note that $\|\mathbf{Q}\times\mathbf{x}\|/\|\mathbf{Q}\|$ is equal to $\|\mathbf{x}\|$ times the sine of the angle that $\mathbf{x}$ forms with the rotation axis.) The radius of the rotation arc at $\mathbf{x}$ and the tangent to the rotation arc at $\mathbf{x}$, together with the bisector of the rotation angle form a right-angle triangle, and, on this triangle, $\|\mathbf{Q}\| = \tan(\theta/2)$ is equal to the length of the opposite side (the tangent) divided by the length of the adjacent side (the radius). Therefore, the length of the tangent, from $\mathbf{x}$ to the bisector, is $\|\mathbf{Q}\times\mathbf{x}\|$.

Proposition (b) is based on (a) and simple vector addition.

To justify proposition (c), start with (b), displace the origin until all the involved vectors lie on the plane perpendicular to the rotation axis, and then normalize $(1+\mathbf{Q}\times)\mathbf{a}$ so that it has unit norm and its head lies on the (unit) rotation arc.

## 4 Derivation of Rotation Formulas

We have already seen in proposition (a) that $\mathbf{Q} \times \mathbf{x}$ connects $\mathbf{x}$ to the intersection point between the tangent to the rotation arc and the bisector of the rotation angle. A very similar reasoning, illustrated on Fig. 2, leads to the conclusion that $-\mathbf{Q} \times \mathbf{R}\mathbf{x}$ connects $\mathbf{R}\mathbf{x}$ to the same point. Therefore,

$$(1 + \mathbf{Q} \times)\mathbf{x} = (1 - \mathbf{Q} \times)\mathbf{R}\mathbf{x} \tag{8}$$

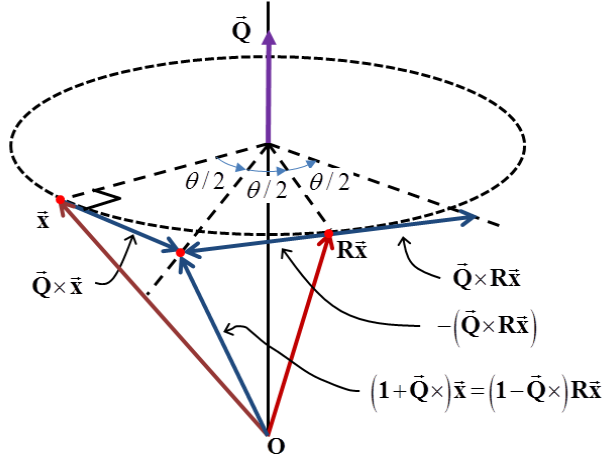

Fig. 2. Geometrical derivation of $(1 + \mathbf{Q} \times)\mathbf{x} = (1 - \mathbf{Q} \times)\mathbf{R}\mathbf{x}$

We only need to assume that $(1 - \mathbf{Q} \times)$ is invertible and left-multiply by its inverse to obtain Cayley's formula for three-dimensional rotations (5). In fact, the inverse matrix exists and can be written explicitly as

$$(1 - \mathbf{Q} \times)^{-1} = 1 + \frac{(\mathbf{Q} \times) + (\mathbf{Q} \times)^2}{1 + \mathbf{Q} \cdot \mathbf{Q}}. \tag{9}$$

This expression can be verified directly, multiplying by $(1 - \mathbf{Q} \times)$ and checking that the result is the unit matrix. During the check one has to apply the properties of the cross product, and, in particular, the identity $(\mathbf{Q} \times)^3 = -(\mathbf{Q} \cdot \mathbf{Q})(\mathbf{Q} \times)$. Inserting Eq. (9) into Eq. (5) and operating (the same identity is needed again) we obtain Eq. (2), which is, in this way, derived from the geometrical interpretation of Rodrigues' vector, without using Euler-Rodrigues' formula. As a matter of fact, since Eq. (1) can be easily recovered from Eq. (2) using the definition of $\mathbf{Q}$, this constitutes a new derivation of Euler-Rodrigues formula, different from existing ones.

The reciprocal of Cayley's rotation formula can be easily derived as well. The first step is the derivation of the expression $(\mathbf{Q} \times)(\mathbf{R} + 1) = (\mathbf{R} - 1)$, illustrated on Fig. 3. The requisite that -1 is not an eigenvalue of $\mathbf{R}$ guarantees that $(\mathbf{R} + 1)$ is invertible. Right-multiplying by its inverse we obtain the desired result (6).

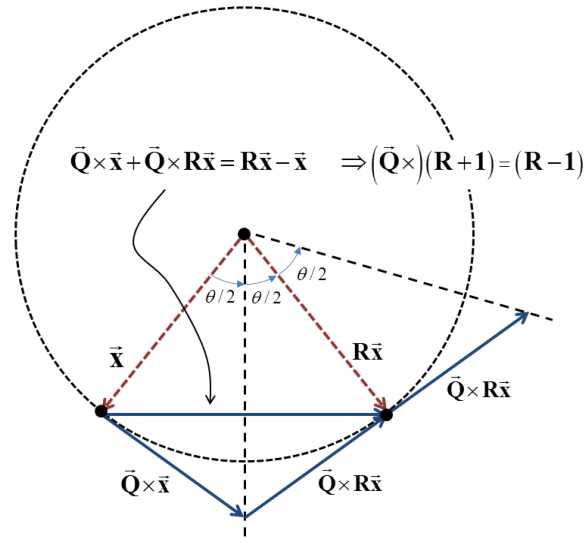

*Fig. 3. Geometrical derivation of $(\mathbf{Q}\times)(\mathbf{R}+\mathbf{1})=(\mathbf{R}-\mathbf{1})$.
(The rotation plane is viewed from above.)*

## 5 Geometrical Derivation of the Rotation Composition Formula

Two consecutive rotations result in other rotation. Formulas that relate the parameters of the resulting rotation to those of the component rotations can be derived in different ways: vector methods [4,16], matrix methods [17], quaternion methods [18], SU(2) matrix methods [19], or a combination of vector methods and spherical trigonometry [20].

In our opinion, the approach that provides the simplest and clearest geometrical visualization of rotation composition, and of the role that half-angles play on it, is a geometrical construction that finds the resultant of two consecutive rotations of a sphere about different diameters, due to Donkin [21]. However, because it is purely geometrical, it does not provide an algebraic formula for the parameters of the resultant rotation. In this section we will review this approach and combine it with the geometrical interpretation of Rodrigues' vector to derive the rotation composition formula (7) in a new, simple and geometrically meaningful way.

Consider a unit sphere, i.e., a sphere of radius 1. A "great circle" is the circle that you get when you intersect a sphere with a plane that contains the center of the sphere. A great arc is any part of the circumference of a great circle.

A rotation about a diameter of a sphere can be specified by a great arc of the great circle perpendicular to this diameter. Thus, being A and B two points on this great circle, a rotation that brings a point from A to B can be represented by the great arc AB. The order of the letters determines the sense of rotation, but the precise position of the arc on its great circle is irrelevant, because it is the length of the great arc what is important. In fact, on the unit sphere, and measuring angles in radians, the length of a great arc is equal to its angle.

With these provisions, Donkin proved the following theorem [21], usually called Donkin's theorem. Let ABC be any spherical triangle. Then, twice the rotation AB, followed by twice the rotation BC, produces the same effect as twice the rotation AC.

Let us reproduce with minor changes Penrose's unbeatably compact justification of Donkin's theorem [22], illustrated in Fig. 4. The spherical triangle ABC is the one mentioned in the theorem, and the other three congruent triangles are the respective reflections in its three vertices. Suppose that a first rotation takes triangle 1 into triangle 2 and a second rotation takes triangle 2 into triangle 3. Then, the rotation that is the composition of the two takes triangle 1 into triangle 3. For triangles 1, 2 and 3 to be actually symmetric to triangle ABC, congruent, and image of each other by rotation, each of these rotations must be through an angle which is precisely twice the corresponding arc-length on the original spherical triangle ABC.

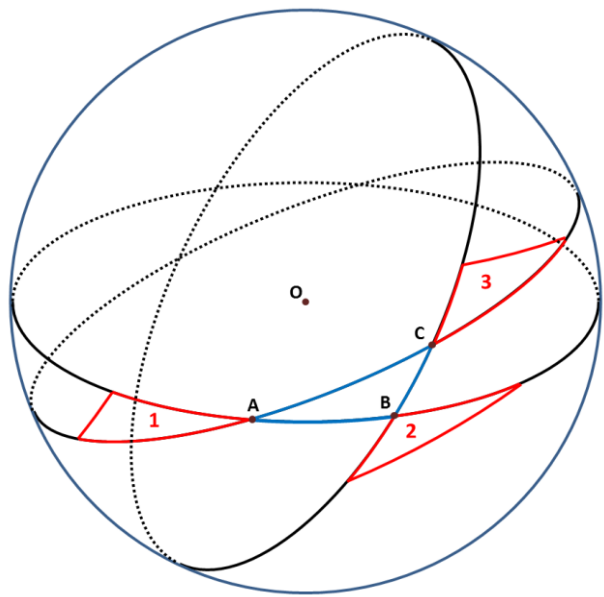

*Fig. 4. Justification of Donkin's Theorem*

Let us compare to translation composition. To compose two translations we represent each one by a vector and then employ vector addition to find the composed translation. This is the "triangle law" illustrated in Fig. 5a. According to Donkin's theorem, rotation composition follows an analogous "spherical triangle law", with one important difference: translation distances correspond to the lengths of the triangle sides, while rotation angles correspond to twice the lengths of the spherical triangle sides. This is illustrated in Fig. 5b.

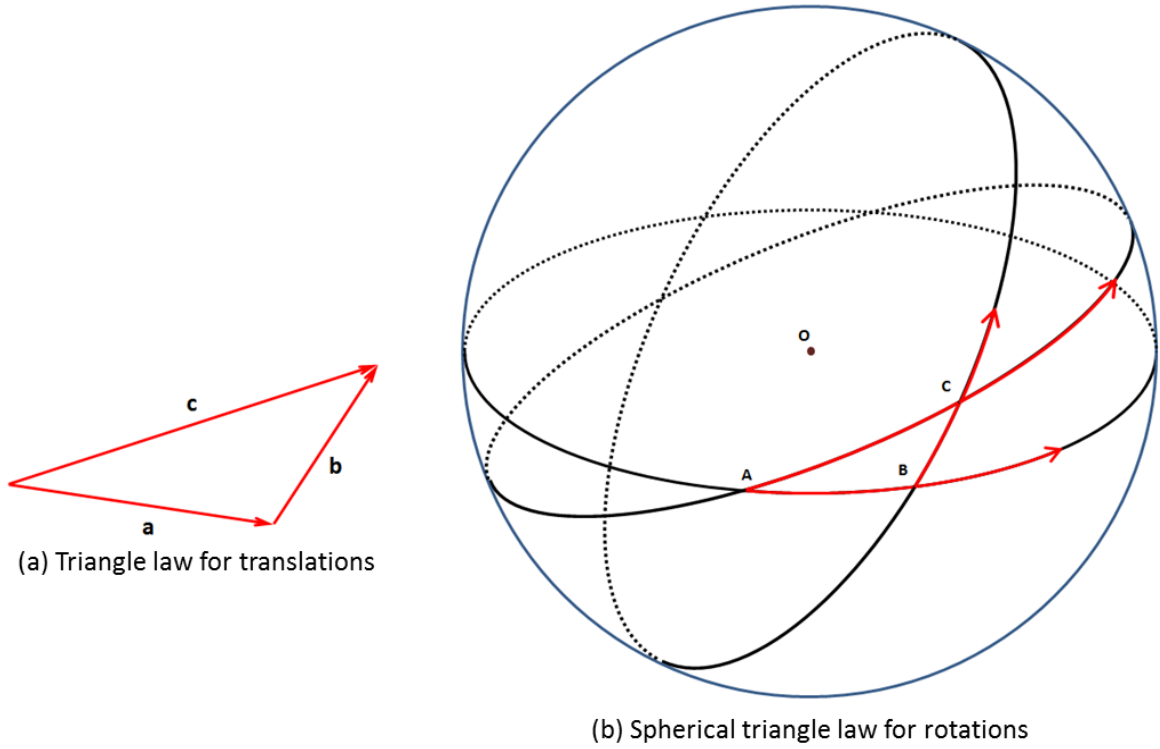

Fig. 5. Rotation composition vs. translation composition

Let **A**, **B** and **C** be the position vectors of the vertices of the spherical triangle on Fig. 5b, and $\mathbf{Q}_1$, $\mathbf{Q}_2$ and $\mathbf{Q}_3$ the Rodrigues' vectors of the first, second and resultant rotations respectively. Our geometrical interpretation of the Rodrigues' vector, and in particular our proposition (c), provides an operator that converts **A** into **B** (based on $\mathbf{Q}_1$), **B** into **C** (based on $\mathbf{Q}_2$) and **A** into **C** (based on $\mathbf{Q}_3$). That is,

$$\frac{(\mathbf{1}+\mathbf{Q}_1\times)\mathbf{A}}{\|(\mathbf{1}+\mathbf{Q}_1\times)\mathbf{A}\|}=\mathbf{B}, \tag{10}$$

$$\frac{(\mathbf{1}+\mathbf{Q}_2\times)\mathbf{B}}{\|(\mathbf{1}+\mathbf{Q}_2\times)\mathbf{B}\|}=\mathbf{C}, \tag{11}$$

$$\frac{(\mathbf{1}+\mathbf{Q}_3\times)\mathbf{A}}{\|(\mathbf{1}+\mathbf{Q}_3\times)\mathbf{A}\|}=\mathbf{C}. \tag{12}$$

These three expressions translate the "spherical triangle law" to vector algebra. We only need a bit of calculation to obtain (7) from them.

Inserting Eq. (10) into Eq. (11) and operating

$$\frac{\left[\mathbf{1}+(\mathbf{Q}_1\times)+(\mathbf{Q}_2\times)+(\mathbf{Q}_2\times)(\mathbf{Q}_1\times)\right]\mathbf{A}}{\|\left[\mathbf{1}+(\mathbf{Q}_1\times)+(\mathbf{Q}_2\times)+(\mathbf{Q}_2\times)(\mathbf{Q}_1\times)\right]\mathbf{A}\|}=\mathbf{C}. \tag{13}$$

Comparing this result to Eq. (12) we see that the numerators and denominators on the left hand side of Eq. (12) and Eq. (13) must be proportional, with the same proportionality constant, let us say $\lambda$, that is,

$$\lambda(\mathbf{1}+\mathbf{Q}_3\times)\mathbf{A}=\left[\mathbf{1}+(\mathbf{Q}_1\times)+(\mathbf{Q}_2\times)+(\mathbf{Q}_2\times)(\mathbf{Q}_1\times)\right]\mathbf{A}. \tag{14}$$

Using the general identity $(\mathbf{Q}_2 \times)(\mathbf{Q}_1 \times) = \left[(\mathbf{Q}_2 \times \mathbf{Q}_1) \times\right] + \mathbf{Q}_2 \mathbf{Q}_1^T - (\mathbf{Q}_2 \cdot \mathbf{Q}_1)\mathbf{1}$, operating, and exchanging terms from one side to the other appropriately,

$$\lambda(\mathbf{Q}_3 \times \mathbf{A}) + (\mathbf{Q}_2 \cdot \mathbf{Q}_1 - 1 + \lambda)\mathbf{A} - \mathbf{Q}_2(\mathbf{Q}_1 \cdot \mathbf{A}) = $$
$$= (\mathbf{Q}_1 \times \mathbf{A}) + (\mathbf{Q}_2 \times \mathbf{A}) + (\mathbf{Q}_2 \times \mathbf{Q}_1) \times \mathbf{A}. \quad (15)$$

Taking into account that $\mathbf{A}$ is perpendicular to $\mathbf{Q}_1$ we can omit the term $\mathbf{Q}_2(\mathbf{Q}_1 \cdot \mathbf{A})$. All the remaining summands are perpendicular to $\mathbf{A}$ (because they are a cross product of $\mathbf{A}$ with other vector) except $(\mathbf{Q}_2 \cdot \mathbf{Q}_1 - 1 + \lambda)\mathbf{A}$, that is obviously parallel to $\mathbf{A}$. The expression can only be true if this last summand is null, that is, if $\lambda = 1 - \mathbf{Q}_2 \cdot \mathbf{Q}_1$. Inserting this value of $\lambda$ into Eq. (15) we get

$$(1 - \mathbf{Q}_2 \cdot \mathbf{Q}_1)(\mathbf{Q}_3 \times \mathbf{A}) = (\mathbf{Q}_1 \times \mathbf{A}) + (\mathbf{Q}_2 \times \mathbf{A}) + (\mathbf{Q}_2 \times \mathbf{Q}_1) \times \mathbf{A}. \quad (16)$$

Hence

$$\mathbf{Q}_3 \times \mathbf{A} = \frac{\mathbf{Q}_1 + \mathbf{Q}_2 + \mathbf{Q}_2 \times \mathbf{Q}_1}{1 - \mathbf{Q}_2 \cdot \mathbf{Q}_1} \times \mathbf{A}. \quad (17)$$

And, if this is true for any vector $\mathbf{A}$ (or at least for any unit vector $\mathbf{A}$ perpendicular to $\mathbf{Q}_1$, which, strictly speaking, is what we have proved) we get Eq. (7), as we wanted to prove.

## 6 Application to Infinitesimal Rotations

For infinitesimal rotations we can neglect second order terms in $\mathbf{Q}$, and rotation formula (2) reduces to

$$\mathbf{R} = \mathbf{1} + 2(\mathbf{Q} \times). \quad (18)$$

To make it clear that we are dealing with infinitesimal rotations, we can write (18) as

$$\mathbf{x} + d\mathbf{x} = (\mathbf{1} + 2\mathbf{Q} \times)\mathbf{x}. \quad (19)$$

Formulas (18) and (19) are very easy to illustrate with our geometrical interpretation of Rodrigues' vector (proposition b and Fig. 1b). In the limit of infinitesimal rotations, the rotation arc and its tangent at point $\mathbf{x}$ coincide, and, while the operator $(\mathbf{1} + \mathbf{Q} \times)$ takes $\mathbf{x}$ halfway to its rotated image, the operator $(\mathbf{1} + 2\mathbf{Q} \times)$ takes $\mathbf{x}$ all the way to $\mathbf{R}\mathbf{x} = \mathbf{x} + d\mathbf{x}$.

For the composition of infinitesimal rotations, neglecting second order terms, the composition law (7) reduces to

$$\mathbf{Q}_3 = \mathbf{Q}_1 + \mathbf{Q}_2. \quad (20)$$

The composition of infinitesimal rotations can thus be performed by adding Rodrigues' vectors. It is obvious from Eq. (20) that infinitesimal rotations commute (although finite ones do not.)

# 7   Application to Rigid Body Kinematics

When a rigid body moves with one point held fixed, it is well known that there exists at each instant of time an instantaneous axis of rotation through the fixed point and an instantaneous angular velocity vector **ω** "along" that axis, such that, taking the fixed point as origin of coordinates, the velocity of any point **x** of the body is given by [23]

$$\frac{d\mathbf{x}}{dt} = \boldsymbol{\omega} \times \mathbf{x}. \qquad (21)$$

If the motion occurs without there being a fixed point in the body, a reference point can be arbitrarily chosen in the body and the movement divided in two components: the translational movement of the reference point with respect to a stationary frame (laboratory frame), and the rotational (or "spinning") movement of the body with respect to a reference frame parallel to that of the laboratory, and whose origin is the reference point. Equation (21) is still valid in the frame in which the reference point is fixed, but needs to be modified, adding the velocity of the reference point, in the laboratory frame.

Comparing equations (19) and (21) it is clear that

$$\mathbf{Q} = \frac{1}{2}\boldsymbol{\omega} dt. \qquad (22)$$

The instantaneous movement of a rigid body spinning around the origin with angular velocity **ω**, during the time interval $dt$, can thus be described as an infinitesimal rotation with Rodrigues' vector $\boldsymbol{\omega} dt/2$.

An important result follows immediately from equations (20) and (22): angular velocities add like vectors in the composition of spinning movements.